\documentclass{amsart}

\usepackage[mathscr]{eucal}
\usepackage{amssymb}
\usepackage{amsthm}
\usepackage{amsmath}

\numberwithin{equation}{section}

\usepackage[all]{xy}
\CompileMatrices

\newdir{ >}{{}*!/-10pt/\dir{>}}

\newtheorem{Thm}{Theorem}[section]
\newtheorem*{Thm*}{Theorem}
\newtheorem{Prop}[Thm]{Proposition}

\newtheorem{Lem}[Thm]{Lemma}
\newtheorem{Cor}[Thm]{Corollary}

\theoremstyle{definition}
\newtheorem{Def}[Thm]{Definition}

\newtheorem{Not}[Thm]{Notation}

\theoremstyle{remark}
\newtheorem{Rem}[Thm]{Remark}

\newtheorem{emptythm}[Thm]{}

\newcommand{\Disp}{\displaystyle}

\newcommand{\bbZ}{\mathbb{Z}}

\newcommand{\calO}{\mathcal{O}}

\newcommand{\isoto}{\buildrel \sim\over\to}
\newcommand{\too}{\mathop{\longrightarrow}}
\newcommand{\isotoo}{\buildrel \sim\over\too}

\newcommand{\cat}[1]{\mathscr{#1}}

\DeclareMathOperator{\Img}{Im}

\DeclareMathOperator{\Spec}{Spec}
\DeclareMathOperator{\Proj}{Proj}
\newcommand{\cone}{{\operatorname{cone}}}

\newcommand{\End}{{\operatorname{End}}}

\newcommand{\inv}{^{-1}}

\newcommand{\pullback}{\ar@{}[rd]|{{\displaystyle\lrcorner}}\ar@{}[rd]|{{\dot{}\,\ }}}
\newcommand{\smallbullet}{{\scriptscriptstyle\bullet}}

\newcommand{\qquadtext}[1]{\qquad\text{#1}\qquad}
\newcommand{\equalby}[2]{\stackrel{#2}{#1}}
\newcommand{\equalbydef}{\equalby{=}{\text{def.}}}

\newcommand{\Mid}{\,\big|\,}
\newcommand{\SET}[2]{\big\{\,#1\Mid#2\,\big\}}
\newcommand{\suchthat}{\text{\,s.t.\,}}
\newcommand{\oursetminus}{\!\smallsetminus\!}
\newcommand{\ourfrac}[2]{\genfrac{}{}{0pt}{}{\scriptstyle #1}{\scriptstyle #2}}

\newcommand{\adhoc}{{\sl ad hoc}}

\newcommand{\apriori}{{\sl a priori}}

\newcommand{\eg}{{\sl e.g.}}

\newcommand{\ie}{{\sl i.e.}\ }
\newcommand{\ipfa}{{\sl ipso facto}}

\newcommand{\mutmut}{{\sl mutatis mutandis}}

\newcommand{\tcat}[1]{(\cat{#1},\otimes,\unit)}
\newcommand{\tcatu}[1]{(\cat{#1},\otimes,\unit_{\cat #1})}
\newcommand{\tenscat}{tensor triangulated category}
\newcommand{\Tenscat}[1]{\tenscat\ $\tcat{#1}$}
\newcommand{\tenscats}{tensor triangulated categories}
\newcommand{\thsub}{thick subcategory}

\newcommand{\tthsub}{thick $\otimes$-ideal}
\newcommand{\tthsubs}{thick $\otimes$-ideals}

\newcommand{\adh}[1]{\overline{#1}}
\newcommand{\adhpt}[1]{\adh{\{#1\}}}

\newcommand{\coll}[1]{\mathcal{#1}}
\newcommand{\iccat}[1]{\widetilde{\cat #1}}
\newcommand{\ideal}[1]{\langle #1\rangle}
\newcommand{\tideal}[1]{\ideal{#1}}
\newcommand{\KZ}[2]{\cat{#1}_{#2}}

\newcommand{\unit}{{1}}
\newcommand{\potimes}[1]{^{\otimes #1}}

\DeclareMathOperator{\Der}{D}
\DeclareMathOperator{\Komp}{K}
\DeclareMathOperator{\smallb}{b}
\DeclareMathOperator{\smallperf}{perf}
\newcommand{\Kb}{\Komp^{\smallb}}
\newcommand{\Dperf}{\Der^{\smallperf}}
\DeclareMathOperator{\proj}{proj}

\DeclareMathOperator{\stab}{stab}

\DeclareMathOperator{\Spc}{Spc}
\newcommand{\Spccat}[1]{\Spc(\cat #1)}

\DeclareMathOperator{\supp}{supp}
\DeclareMathOperator{\supph}{supph}

\DeclareMathOperator{\UU}{U}
\DeclareMathOperator{\ZZ}{Z}
\newcommand{\UUcat}[1]{\UU(\cat #1)}
\newcommand{\ZZcat}[1]{\ZZ(\cat #1)}
\newcommand{\suppcat}[1]{\supp(\cat #1)}
\newcommand{\supponcat}[1]{\supp_{\cat #1}}
\newcommand{\calOcat}[1]{\calO_{\cat #1}}
\newcommand{\Speccat}[1]{\Spec(\cat #1)}

\newcommand{\xytriangle}[7]{
\xymatrix@C=#7em{#1\ar[r]^-{\Disp #4}&#2\ar[r]^-{\Disp #5}&#3\ar[r]^-{\Disp #6}&T #1}}

\newcommand{\xyTriangle}[8]{
\xymatrix@C=#8em{#1\ar[r]^-{\Disp #5}&#2\ar[r]^-{\Disp #6}&#3\ar[r]^-{\Disp #7}&#4}}

\begin{document}


\title[Spectrum of prime ideals in triangulated categories]
{The spectrum of prime ideals\\
in tensor triangulated categories}
\author{Paul Balmer}
\date{\today}

\address{Paul Balmer, D-Math, ETH Zentrum, 8092 Z\"urich, Switzerland}
\email{balmer@math.ethz.ch}
\urladdr{http://www.math.ethz.ch/$\sim$balmer}

\begin{abstract}
We define the \emph{spectrum} of a tensor triangulated category $\cat{K}$ as the set of so-called \emph{prime ideals}, endowed with a suitable topology. In this very generality, the spectrum is the universal space in which one can define \emph{supports} for objects of $\cat{K}$. This construction is functorial with respect to all tensor triangulated functors. Several elementary properties of schemes hold for such spaces, \eg\ the existence of generic points or some quasi-compactness. Locally trivial morphisms are proved to be nilpotent. We establish in complete generality a classification of thick $\otimes$-ideal subcategories in terms of arbitrary unions of closed subsets with quasi-compact complements (Thomason's theorem for schemes, \mutmut). We also equip this spectrum with a sheaf of rings, turning it into a locally ringed space. We compute examples and show that our spectrum unifies the schemes of algebraic geometry and the support varieties of modular representation theory. 
\end{abstract}

\thanks{Research supported by Swiss National Science Foundation, grant~620-66065.}

\maketitle


\vskip-\baselineskip
\vskip-\baselineskip
\section*{Introduction}
\bigbreak

Several mathematicians brought to light the amazing analogies between \apriori\ distinct theories, by means of the triangulated categories naturally appearing in these different areas and more precisely via the so-called \emph{classification of thick subcategories}. Initiated in homotopy theory, see Devinatz, Hopkins and Smith~\cite{dhs}, this classification was transposed to algebraic geometry by Hopkins, see~\cite{hop}, Neeman~\cite{ne92} and Thomason~\cite{thom}. An analogous classification, in terms of support varieties, has been achieved in modular representation theory by Benson, Carlson and Rickard~\cite{bcr} and extended to finite group schemes by Friedlander and Pevtsova~\cite{fripev}.

The importance of triangulated categories is becoming more and more visible all over mathematics. Forged in homological algebra (Grothendieck-Verdier) and in topology (Puppe), these concepts gradually invaded algebraic geometry and marched towards modular representation theory shortly after. Always hidden in the shadow of every newly born Quillen model structure, which also entered complex geometry (see L\'arusson~\cite{fin}), triangulated categories notoriously appeared at the front line of the motivic battle, where Voevodsky won his titles of glory. Recently coveting more analytic lands, they successfully besieged the $C^*$-stronghold of Kasparov's $KK$-theory, see Thom's thesis~\cite{abt}, following a path opened by Higson, Cuntz and others\,; compare Meyer and Nest~\cite{mene}. So, few are the mathematicians who can still be sure that no triangulated category is floating in their ink-pot.

This being said, it is commonly admitted that a triangular structure alone is too arid to be successfully cultivated, as illustrated in $K$-theory. So for irrigation, we adopt the axiomatic of \emph{tensor} triangulated categories $(K,\otimes,\unit)$, since a bi-exact symmetric monoidal structure $\otimes:\cat{K}\times\cat{K}\to\cat{K}$ usually comes along with the triangulation in essentially all known examples. See Def.\,\ref{tens-def}.

In the most recent triangular colonies, a classification of thick subcategories has not yet been achieved, nor even initiated. Nevertheless, it is obvious from the older examples that these objects deserve a study of their own. Many sources deal with triangulated categories\,: Apart from the original Verdier~\cite{verd}, the reader can find a systematic treatise in Neeman~\cite{ne01} and more advanced material in Hovey, Palmieri and Strickland~\cite{hps}. The latter is strongly inspired by classical stable homotopy theory and deals with tensor triangulated categories admitting infinite coproducts in order to use Brown representability, compare~\cite{ne01}. Such additional assumptions have their cost though, in the constant eye one has to keep open for set-theoretical pitfalls. Our lighter setting seems to avoid this burden, so far. Another source is Rosenberg's preprint series~\cite{ros} on non-commutative geometry. Although its focus seems to be more on abelian rather than triangulated categories, it still contains definitions of spectra for triangulated categories, see~\cite[\S\,12]{ros}. However, no recourse to tensor products is made there and, \ipfa, these definitions differ from ours. Further comparison with other trends in non-commutative geometry seems beyond the scope of this introduction (and far beyond the author's present knowledge). Finally, let us mention the \adhoc\ definition of spectrum that Balmer~\cite{bal} used to reconstruct a scheme from its derived category of perfect complexes. As prudently announced there, the definition of~\cite{bal} needs modification before extension outside algebraic geometry. We believe that the construction given below is the right improvement, as we now explain.

\smallbreak

Here, we basically introduce only one ``new'' concept\,: We call \emph{prime} a proper thick $\otimes$-ideal $\cat{P}\subsetneq\cat{K}$ which contains a product only if it contains one of the factors\,:  $a\otimes b\in\cat{P}\ \Longrightarrow\ a\in\cat{P}$ or $b\in\cat{P}$. Reminding us of the super-well-known algebraic notion, this definition certainly passes the first test in being highly digestible. The second test is conceptual quality. It is nonsense to merely mimick some nice definition from commutative algebra, nor from anywhere else as a matter of fact, and then pray for the best. It would be particularly naive here since tensor triangulated categories are not specific to one single area of mathematics, as explained above. We rather need to establish the intrinsic value of our theory for tensor triangulated categories themselves. In this spirit, we prove that the spectrum of $\cat{K}$ is the best locus in which to construct supports for objects of $\cat{K}$ (see precise statement below). Of course, we could \emph{define} the spectrum by this universal property and \emph{then prove} that it can be constructed via prime ideals. We proceed in the reverse order because we expect every reader to be at ease with the notion of prime ideal. The output of the theory for the various examples is the third and final test. 

\smallbreak

We assume that our category $\cat{K}$ is essentially small. We denote by $\Spccat{K}$ the set (!) of prime \tthsubs\ of $\cat{K}$. For any object $a\in\cat{K}$ we denote by $\supp(a):=\{\cat{P}\in\Spccat{K}\mid a\notin\cat{P}\}\subset\Spccat{K}$ \emph{the support of~$a$}. The \emph{Zariski topology} on $\Spccat{K}$ is the one generated by the following basis of open subsets\,: $\{\ \UU(a):=\Spccat{K}\oursetminus\supp(a)\ \mid\ a\in\cat{K}\ \}$. As announced above, Theorem~\ref{univ-thm} says\,:

\begin{Thm*}[\textbf{Universal property of the spectrum}]
We have 
\begin{enumerate}
\item
$\supp(0)=\varnothing$ and $\supp(\unit)=\Spccat{K}$.
\smallbreak
\item
$\supp(a\oplus b)=\supp(a)\cup\supp(b)$.
\smallbreak
\item
$\supp(Ta)=\supp(a)$, for $T:\cat{K}\to \cat{K}$ the translation (shift, suspension).
\smallbreak
\item
$\supp(a)\subset\supp(b)\,\cup\,\supp(c)$ for any exact triangle $a\!\to b\to c\to Ta$.
\smallbreak
\item
$\supp(a\otimes b)=\supp(a)\cap\supp(b)$.
\end{enumerate}
Moreover, for any pair $(X,\sigma)$, where $X$ is a topological space and $\sigma$ an assignment of closed subsets $\sigma(a)\subset X$ to objects $a\in\cat{K}$, which satisfy (a)-(e) above (namely, $\sigma(0)=\varnothing$, $\sigma(\unit)=X$, etc.), there exists a unique continuous map $f:X\to \Spccat{K}$ such that $\sigma(a)=f\inv(\supp(a))$.
\end{Thm*}

The spectrum $\Spc(-)$ is a contravariant functor, for a tensor triangulated functor $F:\cat{K}\to\cat{L}$ induces a map $\Spc(F):\Spccat{L}\to \Spccat{K}$ via $\cat{Q}\mapsto F\inv(\cat{Q})$. In Proposition~\ref{nonvoid-prop}, we prove that $\Spccat{K}$ is non-empty as soon as $\cat{K}\neq0$ and in Corollary~\ref{cl-pts-cor} that $\Spccat{K}$ even has closed points. Although such details might seem trivial to unacquainted readers, let us remind them that former constructions (\eg\ by means of ``atomic'' subcategories, see~\cite{bal}, or by means of indecomposable objects) did not always produce non-empty spaces and were hardly functorial.

Sections~\ref{prime-sect} and~\ref{general-sect} contain other basic results which hopefully illustrate the internal harmony of the theory and which are used in Section~\ref{class-sect} to prove the classification of thick $\otimes$-ideals. For a scheme $X$, Thomason~\cite[Thm.\,3.15]{thom} classifies \tthsubs\ of $\Dperf(X)$ via subsets $Y\subset X$, by assigning to $Y$ the subcategory $\Dperf_{Y}(X)$ of those objects whose homological support is contained in~$Y$. For an arbitrary tensor triangulated category~$\cat{K}$, since $\Spccat{K}$ is the universal locus for supports, the general version of Thomason's Theorem should involve the following construction\,:
$$
\Spccat{K}\supset\qquad Y\ \longmapsto\ \KZ{K}{Y}:=\{a\in\cat{K}\mid\supp(a)\subset Y\}\qquad \subset\cat{K}\,.
$$
Observe though that a subcategory of the form $\KZ{K}{Y}$ is necessarily \emph{radical}, in the usual sense that $a\potimes{n}\in\KZ{K}{Y}\Rightarrow a\in\KZ{K}{Y}$. This follows from $\supp(a\potimes{n})=\supp(a)$, see property~(e) above. Usually, in examples, \emph{all} \tthsubs\ are radical, see Remark~\ref{rad-rem} and Proposition~\ref{rad-prop}, but we do not see any reason for this to hold in general. Hence the statement of Theorem\,\ref{class-thm}\,:

\begin{Thm*}[\textbf{Classification of thick $\boldsymbol{\otimes}$-ideal subcategories}]
Let $\mathfrak{S}$ be the set of those subsets $Y\subset \Spccat{K}$ which are unions $Y=\bigcup_{i\in I}Y_i$ of closed subsets $Y_i$ with quasi-compact complement $\Spccat{K}\oursetminus Y_i$ for all $i\in I$. Let $\mathfrak{R}$ be the set of radical \tthsubs\ of $\cat{K}$. Then there is an order-preserving bijection $\mathfrak{S}\isoto\mathfrak{R}$ given by $Y\mapsto\KZ{K}{Y}=\{a\in\cat{K}\mid\supp(a)\subset Y\}$ with inverse $\cat{J}\mapsto\suppcat{J}:=\bigcup\limits_{a\in\cat{J}}\supp(a)$.
\end{Thm*}

It is then time for computing examples and this is done in Section~\ref{examples-sect}. Indeed, a kind of converse to the above classification holds, namely, a good classification of thick subcategories yields a description of the spectrum (Theorem~\ref{exogene-thm}). Therefore, from the classifications available in the literature, we can easily describe $\Spccat{K}$ for $\cat{K}=\Dperf(X)$ the derived category of perfect complexes over a topologically noetherian scheme~$X$ and for $\cat{K}=\stab(kG)$ the stable category of finitely generated $kG$-modules modulo projective ones, where $G$ is a finite group and $k$ a field of positive characteristic~$p$ (dividing the order of~$G$). Corollaries~\ref{scheme-cor} and~\ref{modular-cor} contain the identification of $\Spccat{K}$ in those two examples, see statement below.

In Section~\ref{ring-sect}, for a general tensor triangulated category~$\cat{K}$, we equip the space $\Spccat{K}$ with a sheaf of rings denoted $\calOcat{K}$ by means of endormorphisms of the unit. The ringed space
$$
\Speccat{K}:=\big(\Spccat{K}\,,\,\calOcat{K}\big)
$$
is always a \emph{locally ringed space}. For the above examples, Theorem~\ref{examples-thm} gives\,:

\begin{Thm*}
With the above notation and hypotheses, we have isomorphisms of schemes\,:
\begin{enumerate}
\item $\Spec\big(\Dperf(X)\big)\simeq X$.
\smallbreak
\item $\Spec\big(\stab(kG)\big)\simeq \Proj\big(H^\smallbullet(G,k)\big)$.
\end{enumerate}
\end{Thm*}


\bigbreak
\noindent\textbf{Acknowledgments\,: } I thank Eric Friedlander, Bruno Kahn and Ralf Meyer for instructive discussions, Zoran Skoda for the reference to~\cite{ros}, and Ivo Dell'Ambrogio, Stefan Gille and Charles Mitchell for their interest and their comments.


\tableofcontents


\section{Terminology and notation}
\label{basis-sect}
\medbreak


Our triangulated category $\cat{K}$ will be essentially small (or choose a fixed universe in which to work). We denote by $Ta$ the translation of an object $a\in\cat{K}$.

\begin{Def}
\label{tens-def}
Here, a \emph{\tenscat} is a triple $(\cat{K},\otimes,\unit)$ consisting of a triangulated category $\cat{K}$ (see~\cite{verd}), a symmetric monoidal ``tensor'' product $\otimes:\cat{K}\times \cat{K}\to \cat{K}$ which is exact in each variable. The unit is denoted~$\unit$ or $\unit_{\cat{K}}$. A \emph{tensor triangulated functor} $F:\cat{K}\to \cat{L}$ is an exact functor respecting the monoidal structures and sending the unit to the unit, $F(\unit_{\cat{K}})=\unit_{\cat{L}}$, unless otherwise stated. At the present stage of the theory, we do not need the higher axiomatic of May~\cite{may} or Keller and Neeman~\cite{kene}.
\end{Def}

\begin{Def}
\label{subcats-def}
A \emph{thick tensor-ideal} $\cat{A}$ of $\cat{K}$ is a full subcategory containing 0 and such that the following conditions are satisfied\,:
\begin{enumerate}
\item
$\cat{A}$ is \emph{triangulated}\,: for any distinguished triangle $a\to b\to c\to Ta$ in $\cat{K}$ if two out of $a$, $b$ and $c$ belong to $\cat{A}$, then so does the third\,;
\item
$\cat{A}$ is \emph{thick}\,: if an object $a\in\cat{A}$ splits in $\cat{K}$ as $a\simeq b\oplus c$ then both summands $b$ and $c$ belong to $\cat{A}$\,;
\item $\cat{A}$ is a \emph{tensor-ideal}\,: if $a\in\cat{A}$ and $b\in\cat{K}$ then $a\otimes b$ also belongs to $\cat{A}$.
\end{enumerate}
Note that (a) forces $\cat{A}$ to be \emph{replete}, \ie closed under isomorphisms\,: $a\simeq b\in\cat{A}$ $\Rightarrow\ a\in\cat{A}$. Since $\cat{K}$ is essentially small, we only have a \emph{set} of such subcategories.
\end{Def}

\begin{Not}
\label{ideal-not}
The intersection of any family of \tthsubs\ is again a \tthsub. Given a collection $\cat{E}$ of objects in $\cat{K}$, we denote by
$
\tideal{\cat{E}}\ \subset\ \cat{K}
$ 
the smallest \tthsub\ of $\cat{K}$ which contains $\cat{E}$.
\end{Not}


\bigbreak
\section{Prime ideals and Zariski topology}
\label{prime-sect}
\medbreak


The main definition of the paper is the following.

\begin{Def}
\label{spec-def}
We call \emph{prime of $\cat{K}$} a proper \tthsub\ $\cat{P}\subsetneq \cat{K}$ such that
$$
a\otimes b\in\cat{P}\qquad\Longrightarrow\qquad a\in\cat{P}\quad\text{or}\quad b\in\cat{P}\,.
$$
Let the \emph{spectrum} of $\cat{K}$, denoted $\Spccat{K}$, be the set of all primes of $\cat{K}$
$$
\Spccat{K}=\{\,\cat{P}\text{ prime of }\cat{K}\,\}\,.
$$
For any family of objects $\cat{S}\subset \cat{K}$ we denote by $\ZZcat{S}$ the following subset of $\Spccat{K}$\,:
\begin{equation}
\label{ZZ-eq}
\ZZcat{S}=\SET{\cat{P}\in\Spccat{K}}{\cat{S}\cap\cat{P}=\varnothing}\,.
\end{equation}
We clearly have $\bigcap_{j\in J}\ZZcat{S_j}=\ZZ(\bigcup_{j\in J}\cat{S}_j)$ and $\ZZcat{S_1}\cup\ZZcat{S_2}=\ZZ(\cat{S}_1\oplus \cat{S}_2)$ where $\cat{S}_1\oplus \cat{S}_2:=\{a_1\oplus a_2\mid a_i\in \cat{S}_i \text{ for }i=1,2\}$. Since $\ZZcat{K}=\varnothing$ and $\ZZ(\varnothing)=\Spccat{K}$, the collection $\SET{\ZZcat{S}\subset\Spccat{K}}{\cat{S}\subset \cat{K}}$ defines the closed subsets of a topology on $\Spccat{K}$, called the \emph{Zariski topology}. We denote the open complement of $\ZZcat{S}$ by
\begin{equation}
\label{UU-eq}
\UU(\cat{S}):=\Spccat{K}\oursetminus\ZZcat{S}=\SET{\cat{P}\in\Spccat{K}}{\cat{S}\cap\cat{P}\neq\varnothing}\,.
\end{equation}
For any object $a\in\cat{K}$, denote by $\supp(a)$ the following closed subset of~$\Spccat{K}$\,:
\begin{equation}
\label{supp-eq}
\supp(a):=\ZZ(\{a\})=\big\{\,\cat{P}\in\Spccat{K}\ \big|\ a\notin\cat{P}\ \big\}
\end{equation}
which we call the \emph{support} of the object~$a\in\cat{K}$.
\end{Def}

A collection of objects $\cat{S}\subset\cat{K}$ is called \emph{(tensor) multiplicative} if $\unit\in\cat{S}$ and if $a_1\,,\,a_2\in \cat{S}\ \Rightarrow a_1\otimes a_2\in \cat{S}$.

\begin{Lem}
\label{nonvoid-lem}
Let $\cat{K}$ be a non-zero \tenscat.
Let $\cat{J}\subset\cat{K}$ be a \tthsub\ and $\cat{S}\subset\cat{K}$ a $\otimes$-multiplicative family of objects such that $\cat{S}\cap\cat{J}=\varnothing$. Then there exists a prime ideal $\cat{P}\in\Spccat{K}$ such that $\cat{J}\subset\cat{P}$ and $\cat{P}\cap \cat{S}=\varnothing$.
\end{Lem}

\begin{proof}
Consider the collection $\coll{F}$ of those \tthsubs\ $\cat{A}\subset\cat{K}$ satisfying\,:
\begin{enumerate}
\item[(1)] $\cat{A}\cap \cat{S}=\varnothing$\,;\kern5em (2)\ \ $\cat{J}\subset\cat{A}$\,;
\smallbreak
\item[(3)] if $c\in \cat{S}$ and $a\in\cat{K}$ are such that $a\otimes c\in\cat{A}$ then $a\in\cat{A}$.
\end{enumerate}
Let $\cat{A}_0:=\{a\in\cat{K}\mid\exists\,c\in \cat{S}\text{ with }a\otimes c\in\cat{J}\}$. One checks directly that $\cat{A}_0$ is a \tthsub\ satisfying properties (1)-(3), hence $\coll{F}$ is non-empty. By Zorn, there exists an element $\cat{P}\in\coll{F}$ maximal for inclusion, which we claim to be prime. Indeed, assume that $a\otimes b\in\cat{P}$ and that $b\notin\cat{P}$ and let us see that $a\in\cat{P}$. Consider
$$
\cat{A}_1:=\{d\in\cat{K}\mid a\otimes d\in\cat{P}\}\,.
$$
One checks easily that $\cat{A}_1$ is a \tthsub, which contains~$\cat{P}$ properly since $b\in \cat{A}_1\oursetminus\cat{P}$. By maximality of $\cat{P}$ in $\coll{F}$, our subcategory $\cat{A}_1$ does \emph{not} belong to~$\coll{F}$. Since $\cat{A}_1$ clearly satisfies properties (2) and (3), it cannot satisfy property~(1), \ie there is an object $d\in \cat{S}$ with $a\otimes d\in \cat{P}$. Now, using property~(3) for $\cat{P}$, we deduce $a\in\cat{P}$.

Note that in case $\cat{S}=\{\unit\}$, condition~(3) is void and $\cat{P}$ merely is a maximal proper ideal containing $\cat{J}$. This proves~(b) in the statement below.
\end{proof}

\begin{Prop}
\label{nonvoid-prop}
Let $\cat{K}$ be a non-zero \tenscat.
\begin{enumerate}
\item
Let $\cat{S}$ be a $\otimes$-multiplicative collection of objects which does not contain zero. Then there exists a prime ideal $\cat{P}\in\Spccat{K}$ such that $\cat{P}\cap \cat{S}=\varnothing$.
\item
Let $\cat{J}\subsetneq\cat{K}$ is a proper \tthsub. Then there exists a maximal proper \tthsub\ $\cat{M}\subsetneq\cat{K}$ which contains $\cat{J}$.
\item
Maximal proper \tthsubs\ are prime.
\item
The spectrum of $\cat{K}$ is not empty\,: $\Spccat{K}\neq\varnothing$.
\end{enumerate}
\end{Prop}

\begin{proof}
(a) is the case $\cat{J}=0$ of Lemma~\ref{nonvoid-lem}. As mentioned in the above proof, (b) is settled. For~(c), apply Lemma~\ref{nonvoid-lem} to $\cat{J}$ maximal and to $\cat{S}=\{\unit\}$\,: there exists a prime containing $\cat{J}$, hence equal to it. Finally, (d) follows from (a), for instance.
\end{proof}

\begin{Cor}
\label{obj-nil-cor}
An object $a\in\cat{K}$ belongs to all primes, $a\in\bigcap_{\cat{P}\in\Spccat{K}}\cat{P}$, \ie $\UU(a)=\Spccat{K}$, \ie $\supp(a)=\varnothing$, if and only if it is $\otimes$-nilpotent, \ie there exists an $n\geq 1$ such that $a\potimes{n}=0$.
\end{Cor}

\begin{proof}
If $a\potimes{n}=0\in\cat{P}$ then $a\in\cat{P}$ for any prime $\cat{P}$. Conversely, if the object $a$ is not nilpotent then $0\notin\cat{S}:=\{a\potimes{n}\mid n\geq0\}$ and we conclude by Prop.\,\ref{nonvoid-prop}\,(a).
\end{proof}

\begin{Cor}
\label{obj-gen-cor}
An object $a\in\cat{K}$ belongs to no prime, \ie $\UU(a)=\varnothing$, \ie $\supp(a)=\Spccat{K}$, if and only if it generates~$\cat{K}$ as a \tthsub, \ie $\tideal{a}=\cat{K}$.
\end{Cor}

\begin{proof}
If $\tideal{a}=\cat{K}$ then $a$ belongs to no proper \tthsub. Conversely, if $\tideal{a}\subsetneq\cat{K}$ is proper, there exists by Prop.\,\ref{nonvoid-prop}\,(b)-(c) a prime $\cat{M}\in\UU(a)$.
\end{proof}

\begin{Lem}
\label{UU-lem}
The assignment $a\,\mapsto\,\UU(a)\equalbydef\{\cat{P}\mid a\in\cat{P}\}$, from objects of $\cat{K}$ to open subsets of $\Spccat{K}$, satisfies the following properties\,:
\begin{enumerate}
\item
$\UU(0)=\Spccat{K}$ and $\UU(\unit)=\varnothing$.
\smallbreak
\item
$\UU(a\oplus b)=\UU(a)\cap\UU(b)$.
\smallbreak
\item
$\UU(Ta)=\UU(a)$.
\smallbreak
\item
$\UU(a)\supset\UU(b)\,\cap\,\UU(c)$ for any exact triangle $a\!\to b\to c\to Ta$.
\smallbreak
\item
$\UU(a\otimes b)=\UU(a)\cup\UU(b)$.
\end{enumerate}
``Dual'' properties hold for the closed complements (see Def.\,\ref{SD-def} and Thm.\,\ref{univ-thm}).
\end{Lem}

\begin{proof}
The properties of a proper \tthsub, see Def.\,\ref{subcats-def}, yield properties (a) to (d) as well as $\UU(a)\cup\UU(b)\subset\UU(a\otimes b)$ in~(e). The other inclusion in~(e) expresses the fact of being prime, see Def.\,\ref{spec-def}.
\end{proof}

\begin{Rem}
\label{basis-rem}
Since for any $\cat{S}\subset\cat{K}$, we have $\UUcat{S}=\bigcup_{a\in \cat{S}}\UU(a)$, it follows from Lemma~\ref{UU-lem}\,(b) that $\{\UU(a)\,|\,a\in\cat{K}\,\}$ is a \emph{basis} of the topology on $\Spccat{K}$. Equivalently, their complements $\{\supp(a)\,|\,a\in\cat{K}\,\}$ form a basis of closed subsets.
\end{Rem}

\begin{Prop}
\label{closure-prop}
Let $W\subset\Spccat{K}$ be a subset of the spectrum. Its closure is
$$
\adh{W}=\kern-.5em\bigcap_{\ourfrac{\Disp a\in\cat{K}\ \suchthat}{\Disp W\subset\supp(a)}}\kern-.5em\supp(a)\,.
$$
\end{Prop}

\begin{proof}
Given a basis $\coll{B}$ of closed subsets, the closure of a subset $W$ is the intersection of all those $B\in\coll{B}$ such that $W\subset B$. We conclude by Remark~\ref{basis-rem}.
\end{proof}

\begin{Prop}
\label{adh-pts-prop}
For any point $\cat{P}\in\Spccat{K}$ its closure in $\Spccat{K}$ is%
$$
\adhpt{\cat{P}}=\{\,\cat{Q}\in\Spccat{K}\,|\,\cat{Q}\subset\cat{P}\,\}\,.
$$
In particular, if $\adhpt{\cat{P}_{1}}=\adhpt{\cat{P}_{2}}$ then $\cat{P}_{1}=\cat{P}_{2}$. (The space $\Spccat{K}$ is $T_0$.)
\end{Prop}

\begin{proof}
Let $\cat{S}_0:=\cat{K}\oursetminus\cat{P}$. Clearly $\cat{P}\in\ZZcat{S_0}$ and if $\cat{P}\in\ZZcat{S}$ then $\cat{S}\subset\cat{S}_0$ and hence $\ZZcat{S_0}\subset\ZZcat{S}$. So, $\ZZcat{S_0}$ is the smallest closed subset which contains the point $\cat{P}$, \ie $\adhpt{\cat{P}}=\ZZcat{S_0}=\{\cat{Q}\in\Spccat{K}\mid\cat{Q}\subset\cat{P}\}$. The second assertion is immediate.
\end{proof}

\begin{Rem}
\label{updown-rem}
Proposition~\ref{adh-pts-prop} is the first indication of the ``reversal of inclusions'' with respect to commutative algebra where the Zariski closure of a point of the spectrum consists of all the \emph{bigger} prime ideals. We shall see in Section~\ref{examples-sect} that the natural homeomorphism $f:\Spec(R)\isoto \Spc\big(\Kb(R-\proj)\big)$ is order-reversing, that is, if $\mathfrak{p}\subset\mathfrak{q}$, as subsets of $R$, then $f(\mathfrak{p})\supset f(\mathfrak{q})$, as subcategories of $\cat{K}$. In this logic, the proof of Lemma~\ref{nonvoid-lem}, which is inspired by the algebraic existence of maximal ideals, does indeed construct maximal ideals of $\Spccat{K}$ but these would be minimal in $\Spec(R)$ when applied to $\cat{K}=\Kb(R-\proj)$. Fortunately, prime ideals also exist at the other end, as we now prove in the usual way.
\end{Rem}

\begin{Prop}
\label{cl-pts-prop}
If $\cat{K}$ is non-zero, there exists minimal primes in $\cat{K}$. More precisely, for any prime $\cat{P}\subset\cat{K}$, there exists a minimal prime $\cat{P}'\subset \cat{P}$.
\end{Prop}

\begin{proof}
To apply Zorn's Lemma, it suffices to see that for any non-empty chain $\coll{C}\subset\Spccat{K}$, the \tthsub\ $\cat{Q}':=\bigcap_{\cat{Q}\in\coll{C}}\cat{Q}$ is a prime. Assume that $a_1\notin\cat{Q}'$ and $a_2\notin\cat{Q}'$. Then there exist $\cat{Q}_i\in\coll{C}$ such that $a_i\notin\cat{Q}_i$ for $i=1,2$. Since $\coll{C}$ is a chain for inclusion, let $\cat{Q}_0$ be the smallest of $\cat{Q}_1$ and $\cat{Q}_2$. Then $a_1,a_2\notin\cat{Q}_0$ and hence $a_1\otimes a_2\notin\cat{Q}_0$, thus $a_1\otimes a_2\notin\cat{Q}'\subset\cat{Q}_0$. In short, $a,b\notin\cat{Q}'\Rightarrow a\otimes b\notin\cat{Q}'$.
\end{proof}

\begin{Cor}
\label{cl-pts-cor}
If the space $\Spccat{K}$ is not empty, it admits a closed point. More precisely, any non-empty closed subset contains at least one closed point.
\end{Cor}

\begin{proof}
Let $\varnothing\neq Z\subset\Spccat{K}$ be closed and let $\cat{P}\in Z$. There exists a $\cat{P}'\subset\cat{P}$ minimal by Prop.\,\ref{cl-pts-prop}. By Prop.\,\ref{adh-pts-prop}, we have $\cat{P}'\in\adhpt{\cat{P}}\subset Z$ and $\adhpt{\cat{P}'}=\{\cat{P}'\}$.
\end{proof}

\begin{Lem}
\label{cover-lem}
Let $a\in\cat{K}$ be an object and $\cat{S}\subset\cat{K}$ be a collection of objects. We have $\UU(a)\subset\UUcat{S}$, \ie $\ZZcat{S}\subset\supp(a)$, if and only if there exists $b_1,\ldots,b_n\in\cat{S}$ such that $b_1\otimes\cdots\otimes b_n\in\tideal{a}$.
\end{Lem}

\begin{proof}
Let $\cat{S}'$ be the $\otimes$-multiplicative collection made of finite products of elements of $\cat{S}\cup\{\unit\}$. Primes being prime, we have $\UUcat{S}=\UUcat{S'}$. So, it clearly suffices to prove the following claim\,: \textit{For $\cat{S}'$ $\otimes$-multiplicative, $\UU(a)\subset\UUcat{S'}$ if and only if \,$\cat{S}'\cap\tideal{a}\neq\varnothing$}. Clearly, if $\cat{S}'\cap\tideal{a}\neq\varnothing$ then any prime containing $a$ meets $\cat{S}'$, that is, $\UU(a)\subset\UUcat{S'}$. Conversely, if $\cat{S}'\cap\tideal{a}=\varnothing$, then Lemma~\ref{nonvoid-lem} (with $\cat{J}:=\tideal{a}$) guarantees the existence of a prime in $\UU(a)\oursetminus\UUcat{S'}$.
\end{proof}

\begin{Prop}
\label{compact-prop}
The following hold true.
\begin{enumerate}
\item
For any object $a\in\cat{K}$, the open $\UU(a)=\Spccat{K}\oursetminus\supp(a)$ is quasi-compact. 
\smallbreak
\item
Any quasi-compact open of $\Spccat{K}$ is of the form $\UU(a)$ for some $a\in\cat{K}$.
\end{enumerate}
\end{Prop}

\begin{proof}
Consider an open covering $\{\UU(\cat{S}_i)\Mid i\in I\}$ of $\UU(a)$. Let $\cat{S}:=\bigcup_{i\in I}\cat{S}_i$, so that $\UU(a)\subset\bigcup_{i\in I}\UU(\cat{S}_i)=\UU(\cat{S})$. By Lemma~\ref{cover-lem}, there exists $b_1,\ldots,b_n\in\cat{S}$ such that $b_1\otimes\ldots\otimes b_n\in\tideal{a}$, but those finitely many objects $b_1,\ldots,b_n$ already belong to $\bigcup_{i\in I_0}\cat{S}_i$ for some finite subset of indices $I_0\subset I$, hence $\UU(a)\subset\bigcup_{i\in I_0}\UU(\cat{S}_i)$.

\smallbreak

For~(b), let $U=\UUcat{S}$ be a quasi-compact open for some $\cat{S}\subset\cat{K}$. Then $U=\bigcup_{a\in\cat{S}}\UU(a)$ and by quasi-compactness of $U$ there exists $a_1,\ldots,a_n\in\cat{S}$ such that $U=\UU(a_1)\cup\ldots\cup\UU(a_n)=\UU(a_1\otimes\cdots\otimes a_n)$, as was to be shown. 
\end{proof}

\begin{Cor}
\label{compact-cor}
Suppose that $\UUcat{S}=\Spccat{K}$ for some collection of objects $\cat{S}\subset\cat{K}$. Then there exists $b_1,\ldots,b_n\in \cat{S}$ such that $b_1\otimes\cdots\otimes b_n=0$. In particular, the spectrum $\Spccat{K}$ is quasi-compact.
\end{Cor}

\begin{proof}
Take $a=0$ in Lem.\,\ref{cover-lem} and Prop.\,\ref{compact-prop}\,(a) and use $\UU(0)=\Spccat{K}$.
\end{proof}

\begin{Rem}
\label{noeth-rem}
A topological space is called noetherian if any non-empty family of closed subsets has a minimal element. This is equivalent to all open subsets being quasi-compact. Hence Proposition~\ref{compact-prop} gives for $\Spccat{K}$\,:
\end{Rem}

\begin{Cor}
\label{noeth-cor}
The topological space $\Spccat{K}$ is noetherian if and only if any closed subset of $\Spccat{K}$ is the support of an object of $\cat{K}$.
\qed
\end{Cor}

\begin{Prop}
\label{irred-prop}
Non-empty irreducible closed subsets of $\Spccat{K}$ have a unique generic point. Indeed, for $\varnothing\neq Z\subset\Spccat{K}$ closed, the following are equivalent\,:
\begin{enumerate}
\item[(i)]
$Z$ is irreducible.
\smallbreak
\item[(ii)]
For all $a,b\in\cat{K}$, if $\UU(a\oplus b)\cap Z=\varnothing$ then $\UU(a)\cap Z=\varnothing$ or $\UU(b)\cap Z=\varnothing$.
\smallbreak
\item[(iii)]
$\cat{P}:=\{a\in\cat{K}\mid\UU(a)\cap Z\neq\varnothing\}$ is a prime.
\end{enumerate}
Moreover, when these conditions hold, we have $Z=\adhpt{\cat{P}}$.
\end{Prop}

\begin{proof}
Uniqueness of generic points is already in Proposition~\ref{adh-pts-prop}.

\smallbreak

\noindent\textit{(i)$\Rightarrow$(ii)\,:}
$Z$ irreducible means that for any open subsets $U_1, U_2$ in $\Spccat{K}$, if $Z\cap U_1\cap U_2=\varnothing$ then $Z\cap U_1=\varnothing$ or $Z\cap U_2=\varnothing$. This gives~(ii), since $\UU(a\oplus b)=\UU(a)\cap\UU(b)$.

\smallbreak

\noindent\textit{(ii)$\Rightarrow$(iii)\,:}
The assumption~(ii) gives $a,b\in\cat{P}\ \Rightarrow\ a\oplus b\in\cat{P}$. Using this, we see that if $a,b\in\cat{P}$ and $a\to b\to c\to T(a)$ is a distinguished triangle, then $c\in\tideal{a\oplus b}$ hence $\UU(a\oplus b)\subset\UU(c)$ and since $\UU(a\oplus b)\cap Z\neq\varnothing$, we get $\UU(c)\cap Z\neq\varnothing$, \ie $c\in\cat{P}$. The other conditions for $\cat{P}$ to be a prime \tthsub\ are easy by Lemma~\ref{UU-lem}.

\smallbreak

\noindent\textit{(iii)$\Rightarrow$(i)\,:}
Let us prove that $Z=\adhpt{\cat{P}}$ which proves~(i) and the ``moreover part''. Let $\cat{Q}\in Z$. For $a\in\cat{Q}$, we have $\cat{Q}\in \UU(a)\cap Z\neq\varnothing$, hence $a\in\cat{P}$. We have proved $\cat{Q}\subset\cat{P}$, \ie $\cat{Q}\in\adhpt{\cat{P}}$ (Prop.\,\ref{adh-pts-prop}) for any $\cat{Q}\in Z$. So, we have $Z\subset\adhpt{\cat{P}}$. Conversely, it suffices to prove $\cat{P}\in \adh{Z}=Z$. To see this, let $s\in\cat{K}$ be an object such that $Z\subset\supp(s)$. Then $\UU(s)\cap Z=\varnothing$ which means $s\notin\cat{P}$ or equivalently $\cat{P}\in\supp(s)$. In short, $\cat{P}\in\kern-1em\bigcap\limits_{s\in\cat{K},\, Z\subset\supp(s)}\kern-1em\supp(s)=\adh{Z}$ by Proposition~\ref{closure-prop}.
\end{proof}

\begin{Cor}
\label{irred-cor}
The spectrum $\Spccat{K}$ is irreducible if and only if for any $a,b$ such that $\tideal{a\oplus b}=\cat{K}$ one has $\tideal{a}=\cat{K}$ or $\tideal{b}=\cat{K}$.\end{Cor}

\begin{proof}
From Proposition~\ref{irred-prop} applied to $Z=\Spccat{K}$, using Corollary~\ref{obj-gen-cor}.
\end{proof}

\begin{Rem}
\label{obj-nil-rem}
We have already established in Corollary~\ref{obj-nil-cor} that an object $a\in\cat{K}$ is $\otimes$-nilpotent if and only if it belongs to all primes $\cat{P}$, which is the same as saying that its image vanish in the localization $\cat{K}/\cat{P}$ for every $\cat{P}\in\Spccat{K}$. We now want to describe the analogue property for morphisms. For a morphism $f:a\to b$, the notation $f\potimes{n}$ of course means the $n$-fold product $f\otimes \cdots\otimes f:\ a\potimes{n}\to b\potimes{n}$.
\end{Rem}

\begin{Prop}
\label{nilp-prop}
Let $f:a\to b$ be a morphism in $\cat{K}$. Suppose that $f$ vanishes in $\cat{K}/\cat{P}$ for all $\cat{P}\in\Spccat{K}$. Then there exists an integer $n\geq 1$ such that $f\potimes{n}=0$.
\end{Prop}

\begin{Lem}
\label{loc0-lem}
Let $f:a\to b$ be a morphism in $\cat{K}$ and let $\cat{P}\in\Spccat{K}$ be a prime. The following conditions are equivalent\,:
\begin{enumerate}
\item[(i)]
$f$ maps to zero in $\cat{K}/\cat{P}$.
\smallbreak
\item[(ii)]
There exists an object $c\in\cat{P}$ such that $f$ factors via~$c$\,:
$\smash{\vcenter{\xymatrix@C=2em@R=1em{&c\ar@{-->}[rd]\\a\ar[rr]^-{f}\ar@{-->}[ru]&&b}}}$.
\end{enumerate}
\end{Lem}
\bigbreak

\begin{proof}
As for any Verdier localization, the assumption $f=0$ in $\cat{K}/\cat{P}$ implies the existence of a morphism $s:z\to a$ such that $f\circ s=0$ and such that $\cone(s)\in\cat{P}$. By the weak cokernel property of the cone, this implies the wanted factorization with $c=\cone(s)$. Conversely, if $f$ factors via some objects which maps to zero in $\cat{K}/\cat{P}$ then $f$ maps to zero in~$\cat{K}/\cat{P}$.
\end{proof}

\begin{proof}[Proof of Proposition~\ref{nilp-prop}]
By assumption and by Lemma~\ref{loc0-lem}, there exists for every $\cat{P}\in\Spccat{K}$ an object $c_{\cat{P}}\in\cat{P}$ such that $f$ factors via~$c_{\cat{P}}$. We have an open covering $\Spccat{K}=\bigcup_{\cat{P}\in\Spccat{K}}\UU(c_{\cat{P}})=\UU\big(\{c_{\cat{P}}\mid\cat{P}\in\Spccat{K}\}\big)$, see notation in Eq.~\eqref{UU-eq}. By quasi-compactness of the spectrum, see Cor.\,\ref{compact-cor}, there exist finitely many objects $c_1,\ldots,c_n\in\cat{K}$ such that $c_1\otimes\cdots\otimes c_n=0$ and such that $f$ factors via each $c_i$. Then $f\potimes{n}:a\potimes{n}\to b\potimes{n}$ factors via $c_1\otimes\cdots\otimes c_n=0$.
\end{proof}


\bigbreak
\section{Universality, functoriality, localization and cofinality}
\label{general-sect}
\medbreak


%
\begin{Def}
\label{SD-def}
A \emph{support data} on a \Tenscat{K} is a pair $(X,\sigma)$ where $X$ is a topological space and $\sigma$ is an assignment which associates to any object $a\in\cat{K}$ a \emph{closed} subset $\sigma(a)\subset X$ subject to the following rules\,:
\begin{enumerate}
\smallbreak
\item[(SD\,1)]
\quad
$\sigma(0)=\varnothing$ and $\sigma(\unit)=X$
\smallbreak
\item[(SD\,2)]
\quad
$\sigma(a\oplus b)=\sigma(a) \cup \sigma(b)$
\smallbreak
\item[(SD\,3)]
\quad $\sigma(Ta)=\sigma(a)$
\smallbreak
\item[(SD\,4)]
\quad
$\sigma(a)\subset\sigma(b)\cup\,\sigma(c)$ for any distinguished triangle $a\to b\to c\to Ta$\,.
\smallbreak
\item[(SD\,5)]
\quad $\sigma(a\otimes b)= \sigma(a)\cap\sigma(b)$\,.
\end{enumerate}

\medbreak

A \emph{morphism} $f:(X,\sigma)\to (Y,\tau)$ of support data on the same category $\cat{K}$ is a continuous map $f:X\to Y$ such that $\sigma(a)=f\inv(\tau(a))$ for all objects $a\in\cat{K}$. 
Such a morphism is an isomorphism if and only if $f$ is a homeomorphism. 
\end{Def}

We now give a universal property for the spectrum.

\begin{Thm}
\label{univ-thm}
Let $\tcatu{K}$ be a \tenscat. The spectrum $(\Spccat{K},\supp)$ of Def.\,\ref{spec-def} is the final support data on $\cat{K}$ in the sense of~\ref{SD-def}. In other words, $(\Spccat{K},\supp)$ is a support data and for any support data $(X,\sigma)$ on $\cat{K}$ there exists a unique continuous map $f:X\to \Spccat{K}$ such that $\sigma(a)=f\inv(\supp(a))$ for any object $a\in\cat{K}$. Explicitly, the map $f$ is defined, for all $x\in X$, by
$$
f(x)=\{a\in\cat{K}\,|\,x\notin\sigma(a)\}\,.
$$
\end{Thm}

\begin{Lem}
\label{unique-lem}
Let $X$ be a set and let $f_1,f_2:X\to\Spccat{K}$ be two maps such that $f_1\inv(\supp(a))=f_2\inv(\supp(a))$ for all $a\in\cat{K}$. Then $f_1=f_2$.
\end{Lem}

\begin{proof}
Let $x\in X$. Observe first that for any object $a\in\cat{K}$, we have by assumption the equivalence $f_1(x)\in\supp(a)\ \Longleftrightarrow\ f_2(x)\in\supp(a)$. This implies in turn that the following two closed subsets of $\Spccat{K}$ coincide\,: 
$$
\bigcap_{f_1(x)\in\supp(a)}\supp(a)\ =\bigcap_{f_2(x)\in\supp(a)}\supp(a)\,.
$$
But the left-hand side is nothing but $\adhpt{f_1(x)}$ and the right-hand side is  $\adhpt{f_2(x)}$ by Prop.\,\ref{closure-prop}. So, $\adhpt{f_1(x)}=\adhpt{f_2(x)}$ in $\Spccat{K}$, hence $f_1(x)=f_2(x)$ by Prop.\,\ref{adh-pts-prop}.
\end{proof}

\begin{Lem}
\label{KY-lem}
Let $(X,\sigma)$ be a support data on $\cat{K}$ and $Y\subset X$ any subset. Then the full subcategory of $\cat{K}$ with objects $\{a\in\cat{K}\,|\,\sigma(a)\subset Y\}$ is a \tthsub.
\end{Lem}

\begin{proof}
Def.\,\ref{subcats-def} is immediately verified using (SD\,1-5) of Def.\,\ref{SD-def}.
\end{proof}

\begin{proof}[Proof of Theorem~\ref{univ-thm}]

We have seen in Lemma~\ref{UU-lem} that $(\Spccat{K},\supp)$ is a support data. For the universal property, let $(X,\sigma)$ be a support data on $\cat{K}$. Uniqueness of the morphism $f:X\to \Spccat{K}$ follows from Lemma~\ref{unique-lem} so let us check that the announced map $f(x):=\{\,a\in\cat{K}\ |\ x\notin\sigma(a)\,\}$ is as wanted. Applying Lemma~\ref{KY-lem} to $Y=X\oursetminus\{x\}$ we see that $f(x)$ is a \tthsub. To see that $f(x)$ is a prime of $\cat{K}$, take $a\otimes b\in f(x)$\,; this means $x\notin\sigma(a\otimes b)=\sigma(a)\cap\sigma(b)$ and hence $x\notin\sigma(a)$ or $x\notin\sigma(b)$, that is, $a\in f(x)$ or $b\in f(x)$. By definition, see Eq.~\eqref{supp-eq}, we have $f(x)\in\supp(a)\ \Leftrightarrow\ a\notin f(x) \ \Leftrightarrow\ x\in\sigma(a)$, hence $f\inv(\supp(a))=\sigma(a)$. This also gives continuity by definition of the topology on $\Spccat{K}$, see Rem.\,\ref{basis-rem}.
\end{proof}

\begin{Rem}
We now want to change the tensor triangulated category~$\cat{K}$. When the dependency on $\cat{K}$ has to be made explicit, we shall denote the support of an element $a\in\cat{K}$ by $\supponcat{K}(a):=\supp(a)\subset\Spccat{K}$.
\end{Rem}

\begin{Prop}
\label{functor-prop}
The spectrum is functorial. Indeed, given a $\otimes$-triangulated functor $F:\cat{K}\to \cat{L}$, the map 
\begin{align*}
\Spc F\,:\ \Spccat{L} & \too\Spccat{K}\,.
\\
\cat{Q}\quad &\longmapsto F\inv(\cat{Q})
\end{align*}
is well-defined, continuous and for all objects $a\in\cat{K}$, we have
$$
(\Spc F)\inv\big(\supponcat{K}(a)\big)\,=\,\supponcat{L}(F(a))
$$
in $\Spccat{L}$. This defines a contravariant functor $\Spc(-)$ from the category of essentially small \tenscats\ to the category of topological spaces. So, if $F:\tcatu{K}\to\tcatu{L}$ and $G:\tcatu{L}\to \tcatu{M}$ are $\otimes$-triangulated functors, then $\Spc(G\circ F)=\Spc(F)\circ\Spc(G)
\ :\ \Spccat{M}\too\Spccat{K}$.
\end{Prop}

\begin{proof}
This is immediate and left as a familiarizing exercise.
\end{proof}

\begin{Cor}
\label{2-cat-cor}
Suppose that two $\otimes$-triangulated functors $F_1,F_2:\cat{K}\to\cat{L}$ satisfy the following property\,: for any $a\in\cat{K}$ we have $\tideal{F_1(a)}=\tideal{F_2(a)}$ in $\cat{L}$, using notation~\ref{ideal-not}. Then the induced maps on spectra coincide\,: $\Spc F_1=\Spc F_2$. This holds in particular if $F_1$ and $F_2$ are (objectwise) isomorphic functors.
\end{Cor}

\begin{proof}
For $\cat{Q}\in\Spccat{L}$, $i\in\{1,2\}$, $a\in\cat{K}$, we have $a\in F_i\inv(\cat{Q})\Leftrightarrow\tideal{F_i(a)}\subset\cat{Q}$.
\end{proof}

\begin{Cor}
\label{ess-surj-cor}
Suppose that a $\otimes$-triangulated functor $F:\cat{K}\to\cat{L}$ is essentially surjective (\ie any object of $\cat{L}$ is isomorphic to the image by $F$ of an object of $\cat{K}$). Then $\Spc F:\Spccat{L}\to\Spccat{K}$ is injective.
\end{Cor}

\begin{proof}
Any prime $\cat{Q}\subset\cat{L}$ is replete (see Def.\,\ref{subcats-def}) and so $\tideal{F(F\inv(\cat{Q}))}=\cat{Q}$, by assumption on~$F$. Therefore, $F\inv(\cat{Q}_1)=F\inv(\cat{Q}_2)$ forces $\cat{Q}_1=\cat{Q}_2$.
\end{proof}

\begin{Prop}
\label{image-prop}
Let $F:\cat{K}\to\cat{L}$ be a $\otimes$-triangulated functor. Let $\cat{S}$ be the collection of those objects $a\in\cat{K}$ whose image generate $\cat{L}$ as a \tthsub, \ie $\cat{S}=\{a\in\cat{K}\mid\tideal{F(a)}=\cat{L}\}$. Then, in the notation of Definition~\ref{spec-def}, the closure of the image of $\Spc(F):\Spccat{L}\to\Spccat{K}$ is
$$
\adh{\Img(\Spc F)}=\ZZcat{S}\,.
$$
\end{Prop}

\begin{proof}
Let $a\in\cat{K}$. We have $a\in\cat{S}$ if and only if $\supponcat{L}(F(a))=\Spccat{L}$, by Cor.\,\ref{obj-gen-cor}. Now, since $(\Spc F)\inv(\supponcat{K}(a))=\supponcat{L}(F(a))$ by Prop.\,\ref{functor-prop}, the condition $a\in\cat{S}$ becomes equivalent to $\Img(\Spc F)\subset\supponcat{K}(a)$. Hence by Prop.\,\ref{closure-prop}, we have
$$
\adh{\Img(\Spc F)}
=\kern-1.5em\bigcap_{\ourfrac{a\in\cat{K}\ \suchthat}{\Img(\Spc F)\subset\supponcat{K}(a)}}\kern-2em\supponcat{K}(a)
=\bigcap_{a\in\cat{S}}\supponcat{K}(a)
\equalbydef\ZZcat{S}\,.
$$
\vskip-\baselineskip
\end{proof}

\begin{Rem}
\label{IJK-rem}
Let $\cat{J}\subset\cat{K}$ be a \tthsub\ of a \Tenscat{K}. Consider $q:\cat{K}\too\cat{L}:=\cat{K}/\cat{J}$ the localization functor. Recall from~\cite{verd} that the quotient category $\cat{L}$ has the same objects as $\cat{K}$ and that its morphisms are obtained via calculus of fractions by inverting those morphisms having their cone in~$\cat{J}$. The category $\cat{L}$ inherits a $\otimes$-structure since $\cat{J}$ is a $\otimes$-ideal. We have a so-called \emph{exact sequence} of \tenscats\,:
\begin{equation}
\label{JKL-eq}
0\to\cat{J}\too^{j}\cat{K}\too^{q}\cat{L}\to0\,.
\end{equation}
The functor $q$ is $\otimes$-triangulated (but not $j$, only because $\cat{J}$ does not have a unit).
\end{Rem}

\begin{Prop}
\label{loc-spc-prop}
Let $q:\cat{K}\to\cat{L}=\cat{K}/\cat{J}$ be a localization as in Remark~\ref{IJK-rem}. The map $\Spc(q):\Spccat{L}\to\Spccat{K}$ induces a homeomorphisms between $\Spccat{L}$ and the subspace $\{\cat{P}\in\Spccat{K}\mid\cat{J}\subset\cat{P}\}$ of $\Spccat{K}$ of those primes containing $\cat{J}$.
\end{Prop}

\begin{proof}
Let us denote by $V=\SET{{P}\in\Spccat{K}}{\cat{J}\subset\cat{P}}$. It is clear that for any $\cat{Q}\in\Spccat{L}$, we have $\Spc(q)\,(\cat{Q})=q\inv(\cat{Q})\supset q\inv(0)=\cat{J}$, \ie $\Img(\Spc(q))\subset V$. We already know from Corollary~\ref{ess-surj-cor} that the map $\Spc(q)$ is injective, since $q:\cat{K}\to\cat{L}$ is (essentially) surjective. Conversely, if $\cat{P}\in\Spccat{K}$ contains $\cat{J}$ then $q(\cat{P})$ is a prime of $\cat{L}$ and $q\inv(q(\cat{P}))=\cat{P}$. This is an easy exercise on Verdier localization and is left to the reader. Finally, let $b=q(a)\in\cat{L}$ and let $\cat{P}\in V$. Then $b\in q(\cat{P})$ if and only if $a\in\cat{P}$ and so $\Spc(q)(\ZZ(b))=\ZZ(a)\cap V$ which proves that the continuous bijection $\Spc(q): \Spccat{L}\to V$ is also a closed map. Hence the result.
\end{proof}

\begin{Rem}
\label{ic-tens-rem}
Recall from~\cite{basch} that the idempotent completion $\iccat{K}$ (which exists for any additive category) of a triangulated category~$\cat{K}$ is canonically triangulated in such a way that the functor $\iota:\cat{K}\to\iccat{K}$ is exact. If $\cat{K}$ is a \tenscat, it is easy to turn $\iccat{K}$ into a \tenscat\ as well so that the functor $\iota:\cat{K}\to\iccat{K}$ is $\otimes$-triangulated. We now show that this does not affect the spectrum. We prove this slightly more generally.
\end{Rem}

\begin{Prop}
\label{cofinal-prop}
Let $\cat{L}$ be a \tenscat\ and let $\cat{K}\subset\cat{L}$ be a full tensor triangulated subcategory with the same unit and which is \emph{cofinal}, \ie for any object $a\in\cat{L}$ there exists $a'\in\cat{L}$ such that $a\oplus a'\in\cat{K}$. Then the map $\cat{Q}\mapsto\cat{Q}\cap\cat{K}$ defines a homeomorphism $\Spccat{L}\isoto\Spccat{K}$.
\end{Prop}

\begin{proof}
Replacing $\cat{K}$ by its isomorphic-closure, we can assume that $\cat{K}$ is replete (see Def.\,\ref{subcats-def}). The above map $\Spccat{L}\to\Spccat{K}$ is nothing but $\Spc(\iota)$ where $\iota:\cat{K}\hookrightarrow\cat{L}$ is the inclusion, so it is well-defined and continuous. Recall the well-known fact\,:
\begin{equation}
\label{cofinal-eq}
\text{for any object }a\in\cat{L},\text{ we have }a\oplus T(a)\in\cat{K}\,,
\end{equation}
whose proof we give for the reader's convenience. There exists by assumption an object $a'\in\cat{L}$ such that $a\oplus a'\in\cat{K}$. Let us add to the distinguished triangle $a'\to0\to T(a')\to T(a')$ two other triangles, namely $a\to a\to0\to T(a)$ and $0\to T(a)\to T(a)\to 0$, to obtain the distinguished triangle
$$
\xytriangle{(a\oplus a')}{a\oplus T(a)}{T(a\oplus a')}{}{}{}{3}
$$
which has two entries in $\cat{K}$ and hence the third\,: $a\oplus T(a)\in\cat{K}$. This proof also shows that if $a\oplus a'$ belongs to some triangulated subcategory (\eg\ a prime) of $\cat{K}$ then so does $a\oplus T(a)$. So, given a prime $\cat{P}\in\Spccat{K}$, we have the equality
\begin{equation}
\label{icP-eq}
\{a\in\cat{L}\mid a\oplus T(a)\in\cat{P}\}=\{a\in\cat{L}\mid\exists\, a'\in\cat{L}\suchthat a\oplus a'\in\cat{P}\}=:\widetilde{\cat{P}}\,.
\end{equation}
We claim that $\widetilde{\cat{P}}$ is a prime of $\cat{L}$. It is easy to check that it is a \tthsub. Suppose that $a\otimes b\in\widetilde{\cat{P}}$ and that $a\notin\widetilde{\cat{P}}$. This means that if we let $c:=a\oplus T(a)$ we have $c\in\cat{K}\oursetminus\cat{P}$ and $c\otimes b\simeq (a\otimes b)\oplus T(a\otimes b)\in\cat{P}$. The latter $c\otimes b\in\cat{P}$ implies that $c\otimes (b\oplus T(b))\simeq(c\otimes b)\oplus T(c\otimes b)\in\cat{P}$ and hence that $b\oplus T(b)\in\cat{P}$, since $\cat{P}$ is prime and does not contain~$c$. So, we have proved that $b\in\widetilde{\cat{P}}$, as wanted. Since $\cat{P}$ is thick, it is easy to see that $\widetilde{\cat{P}}\cap\cat{K}=\cat{P}$. So, $\cat{P}\mapsto\widetilde{\cat{P}}$ is a right inverse to $\Spc(\iota)$.

Let $\cat{Q}\in\Spccat{L}$. Then $\cat{Q}=\widetilde{\cat{P}}$ where $\cat{P}:=\cat{Q}\cap\cat{K}\in\Spccat{K}$. The inclusion $\cat{Q}\subset\widetilde{\cat{P}}$ is obvious from Eq.~\eqref{cofinal-eq} above. The other inclusion follows from thickness of~$\cat{Q}$.

So $\Spc(\iota)$ is a continuous bijection with inverse given by $\cat{P}\mapsto\widetilde{\cat{P}}$. Clearly, for any $a\in\cat{L}$ we have $a\in\cat{Q}$ if and only if $a\oplus T(a)\in\cat{Q}$ if and only if $a\oplus T(a)\in\cat{P}$, where $\cat{P}$ and $\cat{Q}$ are corresponding primes, \ie $\cat{P}=\cat{Q}\cap\cat{K}$ and $\cat{Q}=\widetilde{\cat{P}}$. This shows that the image by $\Spc(\iota)$ of the closed subset $\supponcat{L}(a)$ is $\supponcat{K}(a\oplus T(a))$ which is closed. Hence $\Spc(\iota)$ is a closed map (see Rem.\,\ref{basis-rem}).
\end{proof}

\begin{Cor}
\label{ic-cor}
We have a homeomorphism $\Spc(\iota):\Spc(\iccat{K})\isoto\Spccat{K}$ where $\iota:\cat{K}\hookrightarrow\iccat{K}$ is the idempotent completion of~$\cat{K}$, see Remark~\ref{ic-tens-rem}.
\qed
\end{Cor}


\bigbreak
\section{Classification of thick subcategories}
\label{class-sect}
\medbreak


%
\begin{Def}
\label{rad-def}
As usual, the \emph{radical} $\sqrt{\cat{J}}$ of a \tthsub\ $\cat{J}\subset\cat{K}$ is defined to~be
$$
\sqrt{\cat{J}}:=\{a\in\cat{K}\mid\exists\,n\geq1\text{ such that }a\potimes{n}\in\cat{J}\}\,.
$$
A \thsub\ $\cat{J}$ is called \emph{radical} if $\sqrt{\cat{J}}=\cat{J}$.
\end{Def}

\begin{Lem}
\label{sqrtJ-lem}
$\sqrt{J}$ is a \tthsub\ equal to $\bigcap\limits_{\Disp\kern-1em\cat{P}\in\Spccat{K},\,\cat{J}\subset\cat{P}}\kern-1em\cat{P}$\,.
\end{Lem}

\begin{proof}
It suffices to prove the claimed equality $\sqrt{\cat{J}}=\bigcap_{\cat{P}\supset\cat{J}}\cat{P}$ since an intersection of \tthsubs\ still is a \tthsub. Clearly, by definition of primes, $\sqrt{\cat{J}}\subset\cat{P}$ for any prime $\cat{P}$ containing $\cat{J}$. Conversely, let $a\in\cat{K}$ be an object such that $a\in\cat{P}$ for all $\cat{P}\supset\cat{J}$. Consider the $\otimes$-multiplicative $\cat{S}:=\{a\potimes{n}\mid n\geq1\}$. We have to show that $\cat{S}\cap\cat{J}\neq\varnothing$. Indeed, $\cat{S}\cap\cat{J}=\varnothing$ is excluded by Lemma~\ref{nonvoid-lem} which would give us a prime $\cat{P}$ with $\cat{J}\subset\cat{P}$ and $a\notin\cat{P}$, contradicting the choice of the object~$a$.
\end{proof}

\begin{Rem}
\label{rad-rem}
In practice, it is very frequent that \emph{all} \tthsubs\ are radical. Indeed, as soon as an object $a\in\cat{K}$ is dualizable, we have that $a$ is a direct summand of $a\otimes a\otimes D(a)$ where $D(a)$ is the dual of $a$. See details in~\cite[Lem.\,A.2.6]{hps}. In particular, it is very common that $a\in\tideal{a\otimes a}$, in which case we can use\,:
\end{Rem}

\begin{Prop}
\label{rad-prop}
The following conditions are equivalent\,:
\begin{enumerate}
\item[(i)]
Any \tthsub\ of the category $\cat{K}$ is radical.
\smallbreak
\item[(ii)]
We have $a\in\tideal{a\otimes a}$ for all objects $a\in\cat{K}$.
\end{enumerate}
\end{Prop}

\begin{proof}
If any \tthsub\ $\cat{J}$ is radical, then so is $\cat{J}:=\tideal{a\otimes a}$, giving~(ii). Conversely, suppose that~(ii) holds and let $\cat{J}$ be a \tthsub. We have to show that $a\potimes{n}\in\cat{J}$ $\Rightarrow\ a\in\cat{J}$. By induction on~$n$, it suffices to treat the case $n=2$, which is immediate from the assumption\,: $a\in\tideal{a\otimes a}\subset\cat{J}\ \Rightarrow\ a\in\cat{J}$.
\end{proof}

\begin{Not}
\label{suppE-not}
The \emph{support} of a collection of objects $\cat{E}\subset \cat{K}$ is defined to be the union of the supports of its elements\,:
$$
\suppcat{E}=\bigcup_{a\in \cat{E}}\supp(a)\ \subset\,\Spccat{K}\,.
$$ 
Warning\,: the subset $\suppcat{E}\subset\Spccat{K}$ is not the same thing as the closed subset $\ZZcat{E}\subset\Spccat{K}$ of Eq.~\eqref{ZZ-eq}, although both coincide with $\supp(a)$ when $\cat{E}=\{a\}$.
\end{Not}

\goodbreak

\begin{Lem}
\label{suppE-lem}
Let $\cat{\cat{E}}\subset\cat{K}$. Then $\suppcat{E}=\{\cat{P}\in\Spccat{K}\mid\cat{E}\not\subset\cat{P}\}$.
\end{Lem}

\begin{proof}
We have $\cat{P}\in\suppcat{E}$ if and only if there exists an $a\in\cat{E}$ such that $\cat{P}\in\supp(a)$ which means $a\notin\cat{P}$, by definition of the support, see Eq.~\eqref{supp-eq}.
\end{proof}

\begin{Not}
\label{KZ-not}
Given a subset $Y\subset\Spccat{K}$ we define the subcategory \emph{supported on $Y$} to be the full subcategory $\KZ{K}{Y}$ of $\cat{K}$ on the following objects\,:
$$
\KZ{K}{Y}=\{\,a\in\cat{K}\,|\,\supp(a)\subset Y\,\}\ \subset\,\cat{K}\,.
$$
It is a \tthsub\ by Lemma~\ref{KY-lem}.
\end{Not}

\begin{Lem}
\label{KZ-lem}
Let $Y\subset\Spccat{K}$ be a subset. Then $\KZ{K}{Y}=
\bigcap\limits_{\Disp\cat{P}\in \Spccat{K}\oursetminus Y}\cat{P}$.
\end{Lem}

\begin{proof}
This is easy, for instance as follows\,: For an object $a\in\cat{K}$, we have $a\in\KZ{K}{Y}$ if and only if $\supp(a)\subset Y$. Taking complements, $\supp(a)\subset Y$ is in turn equivalent to\, ``\,$\forall\,\cat{P}\notin Y$, $\cat{P}\notin\supp(a)$\,''. Now, by definition of $\supp(a)=\{\cat{P}\in\Spccat{K}\mid a\notin\cat{P}\}$, the latter property of the object~$a$ is equivalent to\,: ``\,$\forall\,\cat{P}\notin Y$, $a\in\cat{P}$\,''.
\end{proof}

\begin{Prop}
\label{rad-class-prop}
Let $\cat{J}\subset\cat{K}$ be a \tthsub. Then we have $\KZ{K}{\suppcat{J}}=\sqrt{\cat{J}}$.
\end{Prop}

\begin{proof}
By Lemma~\ref{KZ-lem}, $\KZ{K}{\suppcat{J}}=\bigcap_{\cat{P}\notin\suppcat{J}}\cat{P}$. By Lemma~\ref{suppE-lem}, $\cat{P}\notin\suppcat{J}$ is equivalent to $\cat{J}\subset\cat{P}$. So, we conclude by Lemma~\ref{sqrtJ-lem}.
\end{proof}

\begin{Thm}
\label{class-thm}
Let $\mathfrak{S}$ be the set of those subsets $Y\subset \Spccat{K}$ of the form $Y=\bigcup_{i\in I}Y_i$ for closed subsets $Y_i$ of $\Spccat{K}$ with $\Spccat{K}\oursetminus Y_i$ quasi-compact for all $i\in I$. Let $\mathfrak{R}$ be the set of radical \tthsubs\ of $\cat{K}$. Then there is an order-preserving bijection $\mathfrak{S}\isoto\mathfrak{R}$ given by
$$
\qquad
Y\longmapsto\KZ{K}{Y}\equalbydef\{a\in\cat{K}\mid\supp(a)\subset Y\}
\qquadtext{(see Notation~\ref{KZ-not})}
$$
whose inverse is 
$$
\qquad
\cat{J}\longmapsto\suppcat{J}\equalbydef\bigcup_{a\in\cat{J}}\supp(a)
\qquad\qquadtext{(see Notation~\ref{suppE-not})\,.}
$$
\end{Thm}

\begin{proof}
Let us see that the maps are well-defined. A \tthsub\ of the form $\KZ{K}{Y}$ is necessarily radical since $\supp(a\potimes{n})=\supp(a)\cap\ldots\cap\,\supp(a)=\supp(a)$ by Lem.\,\ref{UU-lem}\,(e). By Prop.\,\ref{compact-prop}\,(a), $\suppcat{J}$ is a union of closed subsets with quasi-compact complements $\Spccat{K}\oursetminus\supp(a)=\UU(a)$. Both maps are clearly inclusion-preserving. Let us now check that both composite are equal to the identity.

We already know from Prop.\,\ref{rad-class-prop} that $\KZ{K}{\suppcat{J}}=\sqrt{\cat{J}}=\cat{J}$ for $\cat{J}$ radical. We now turn to the other composition\,: $\supp(\KZ{K}{Y})$. Observe that for any subset $Y\subset\Spccat{K}$, we have directly from the definitions that $\supp(\KZ{K}{Y})\subset Y$. Conversely, let $\cat{P}\in Y$ and let us prove that $\cat{P}\in\supp(\KZ{K}{Y})$. By assumption on $Y$, there exists a closed subset $Y_i\subset Y$ such that $\cat{P}\in Y_i$ and $\Spccat{K}\oursetminus Y_i$ is quasi-compact. By Prop.\,\ref{compact-prop}\,(b), there exists an object $a\in\cat{K}$ such that $\Spccat{K}\oursetminus Y_i=\UU(a)$, that is, $Y_i=\supp(a)$. So, we have $\cat{P}\in\supp(a)\subset Y$. This means $\cat{P}\in\bigcup_{a\in\KZ{K}{Y}}\supp(a)=\supp(\KZ{K}{Y})$ as was to be shown.
\end{proof}

\begin{Rem}
\label{re-noeth-rem}
If $\Spccat{K}$ is a noetherian topological space, then one can replace Thomason's condition ``$Y=\bigcup Y_i$ with $\Spccat{K}\oursetminus Y_i$ quasi-compact'' by the simpler condition ``$Y$ specialization closed''. In case of doubt, see Definition~\ref{classif-def} below, or Remark~\ref{noeth-rem} and Corollary~\ref{noeth-cor}.
\end{Rem}

\begin{Rem}
\label{class-K0-rem}
Assume for simplicity that all \tthsubs\ are radical (see Rem.\,\ref{rad-rem} and Prop.\,\ref{rad-prop}). Then, if we want to describe $\otimes$-ideals $\cat{A}$ of $\cat{K}$ which are not necessarily thick, we can do it in two steps. First, consider $\cat{J}:=\{a\in\cat{K}\mid a\oplus T(a)\in\cat{A}\}$ which is the same by~\eqref{cofinal-eq} as $\cat{J}=\{a\in\cat{K}\mid\exists\,b\in\cat{K}\text{ with }a\oplus b\in\cat{A}\}$ and which obviously admits $\cat{A}$ as a cofinal subcategory. This category $\cat{J}$ is a \tthsub\ of $\cat{K}$ and can be classified via its support, as explained above. Finally, $\cat{A}$ can be recovered from $\cat{J}$ via the subgroup it defines in the zeroth $K$-theory group $K_0(\cat{J})$ of $\cat{J}$, as explained in Thomason~\cite[Thm.\,2.1]{thom} (whose ``dense'' is our ``cofinal'').
\end{Rem}


\bigbreak
\section{Examples}
\label{examples-sect}
\medbreak


To avoid repeating several times the same property, we give it a name\,:

\begin{Def}
\label{classif-def}
Recall that a subset $Y\subset X$ of a topological space~$X$ is \emph{specialization closed} if it is a union of closed subsets or equivalently if $y\in Y$ implies $\adhpt{y}\subset Y$.
We say for short that a support data $(X,\sigma)$ on a \tenscat\ $\cat{K}$ (Def.\,\ref{SD-def}) is a \emph{classifying support data} if the following two conditions hold\,:
\begin{enumerate}
\item
The topological space $X$ is noetherian and any non-empty irreducible closed subset $Z\subset X$ has a unique generic point\,: $\exists!\,x\in Z$ with $\adhpt{x}=Z$. 
\smallbreak
\item
We have a bijection $\theta:\{Y\subset X\mid Y\text{ specialization closed}\}\isotoo$\break$\SET{\cat{J}\subset\cat{K}}{\,\cat{J}\text{ radical \tthsub}}$ defined by $Y\mapsto\{a\in\cat{K}\mid\sigma(a)\subset Y\}$, with inverse $\cat{J}\mapsto\sigma(\cat{J}):=\bigcup_{a\in\cat{J}}\sigma(a)$.
\end{enumerate}
\end{Def}

\begin{Thm}
\label{exogene-thm}
Suppose that $(X,\sigma)$ is a classifying support data on~$\cat{K}$. Then the canonical map $f:X\to \Spccat{K}$ of Theorem~\ref{univ-thm} is a homeomorphism.
\end{Thm}

\begin{proof}
Theorem~\ref{univ-thm} tell us that the map $f$ is continuous and satisfies $f\inv(\supp(a))=\sigma(a)$ for all objects $a\in\cat{K}$. We first prove the following\,:

\medbreak

\noindent{\it Claim\,: Any closed subset $Z\subset X$ is of the form $Z=\sigma(a)$ for some object $a\in\cat{K}$.}

\medbreak

Since $\sigma(a_1)\cup\ldots\cup\sigma(a_n)=\sigma(a_1\oplus\ldots\oplus a_n)$ and since the space $X$ is noetherian, it suffices to prove the claim for an irreducible $Z=\adhpt{x}$ for some $x\in X$. Now, $\adhpt{x}=Z=\theta\inv(\theta(Z))=\bigcup_{a\in\theta(Z)}\sigma(a)$ forces the existence of some object $a\in\cat{K}$ such that $x\in\sigma(a)\subset Z$. Hence $\adhpt{x}\subset\sigma(a)\subset Z=\adhpt{x}$ which proves the claim.

\medbreak

For $x\in X$ define $Y(x):=\SET{y\in X}{x\notin\adhpt{y}}$. It is easy to check that $Y(x)$ is specialization closed. Let $a\in \cat{K}$. Let us see that $\sigma(a)\subset Y(x)\Leftrightarrow x\notin\sigma(a)$. Since $x\notin Y(x)$ we have\,: $\sigma(a)\subset Y(x)\Rightarrow x\notin\sigma(a)$. Conversely, since $\sigma(a)$ is (specialization) closed, and since $x\notin\sigma(a)$, we have $x\notin\adhpt{y},\ \forall y\in\sigma(a)$, which exactly means $\sigma(a) \subset Y(x)$ by definition of the latter. So, we have established\,:
\begin{equation}
\label{sigma-1-eq}
\theta(Y(x))\equalbydef\SET{a\in\cat{K}}{\sigma(a)\subset Y(x)}
=\SET{a\in\cat{K}}{x\notin\sigma(a)}\equalbydef f(x)\,.
\end{equation}
In particular, if $f(x_1)=f(x_2)$ then $Y(x_1)=Y(x_2)$ which immediately implies $\adhpt{x_1}=\adhpt{x_2}$ and $x_1=x_2$ since $X$ is $T_0$, see Def.\,\ref{classif-def}\,(a). Hence $f$ is injective.

\smallbreak

Let us check surjectivity of the map~$f$. Let $\cat{P}$ be a prime of $\cat{K}$. By assumption, there exists a specialization closed subset $Y\subset X$ such that $\cat{P}=\theta(Y)$. The complement $X\oursetminus Y$ is non-empty since $\cat{P}\not=\cat{K}$. Let $x,y\in X\oursetminus Y$. By the Claim, there exist objects $a,b\in\cat{K}$ such that $\adhpt{x}=\sigma(a)$ and $\adhpt{y}=\sigma(b)$. Since $x$ and $y$ are outside $Y$, the objects $a$ and $b$ do not belong to $\theta(Y)=\cat{P}$. The latter being a prime, we then have $a\otimes b\notin\cat{P}$. So, $\sigma(a\otimes b)\not\subset Y$, \ie there is a point $z\in X\oursetminus Y$ such that $z\in\sigma(a\otimes b)=\sigma(a)\cap\sigma(b)=\adhpt{x}\cap\adhpt{y}$ and hence $\adhpt{z}\subset\adhpt{x}$ and $\adhpt{z}\subset\adhpt{y}$. In short, we have established that the non-empty family of closed subsets
$$
\coll{F}:=\SET{\adhpt{x}\subset X}{x\in X\oursetminus Y}
$$
has the property that any two elements admit a lower bound for inclusion. On the other hand, since $X$ is noetherian, there exists a minimal element in $\coll{F}$ which is then \emph{the} lower bound for $\coll{F}$ by the above reasoning. This shows that there exists a point $x\in X\oursetminus Y$ such that $X\oursetminus Y\subset\{y\in X\mid x\in\adhpt{y}\}$, the reverse inclusion also holds because $x\notin Y$, which is specialization closed. In other words,
$$
Y=\{y\in X\Mid x\notin\adhpt{y}\}=Y(x)\,.
$$
Therefore $\cat{P}=\theta(Y)=\theta(Y(x))\equalby{=}{\eqref{sigma-1-eq}}f(x)$, which proves the surjectivity of~$f$.

\smallbreak

The relation $f\inv(\supp(a))=\sigma(a)$ now gives $f(\sigma(a))=\supp(a)$. This shows that $f$ is a closed map, since we know from the Claim that any closed subset of $X$ is of the form $\sigma(a)$. Hence the map $f$ is a homeomorphism.
\end{proof}

\begin{Rem}
\label{zut-rem}
Theorem~\ref{exogene-thm} is a converse to Theorem~\ref{class-thm} in the noetherian case. We now want to use it in order to describe $\Spccat{K}$ in two classes of examples.
\end{Rem}

\begin{Not}
\label{scheme-not}
Let $X$ be a (topologically) noetherian scheme and let $\Dperf(X)$ be the derived category of perfect complexes over $X$, with the usual tensor product $\otimes=\otimes^L_{\calO_X}$. For any perfect complex $a\in\Dperf(X)$ we denote by $\supph(a)\subset X$ the \emph{homological support} of $a$, which is the support of the (total) homology of $a$.
\end{Not}

\begin{Thm}[Thomason~{\cite[Thm.\,3.15]{thom}}]
\label{thom-thm}
The pair $(X\,,\,\supph)$ is a classifying support data on $\Dperf(X)$ in the sense of Definition~\ref{classif-def}.
\end{Thm}

\begin{Cor}
\label{scheme-cor}
There is a homeomorphism $f:X\isoto\Spc\big(\Dperf(X)\big)$ with
$$
\qquad 
f(x)=\SET{a\in\Dperf(X)}{\,a_x\simeq 0\text{ in }\Dperf(\calO_{X,x})}\qquad\text{for all }x\in X\,.
$$
Moreover, for any perfect complex $a\in\Dperf(X)$, the closed subset $\supph(a)$ in~$X$ corresponds via~$f$ to the closed subset $\supp(a)$ in $\Spc\big(\Dperf(X)\big)$.
\end{Cor}

\begin{proof}
Theorems~\ref{exogene-thm}, \ref{thom-thm} and~\ref{univ-thm}.
\end{proof}

\begin{Rem}
\label{scheme-rem}
Note that we can not expect any scheme to be recovered as $\Spccat{K}$ from $\cat{K}=\Dperf(X)$, nor from other tensor triangulated categories $\cat{K}$ associated to~$X$, since $\Spccat{K}$ is always quasi-compact (Cor.\,\ref{compact-cor}) and quasi-separated, \ie it admits a basis of quasi-compact open (Rem.\,\ref{basis-rem} and Prop.\,\ref{compact-prop}).
\end{Rem}

\begin{Not}
\label{group-not}
Let $G$ be a finite group or more generally a finite group scheme and let $k$ be a field of characteristic~$p$. We adopt the notations of~\cite{fripev}. Let $H^\smallbullet(G,k)$ be $\oplus_{i\in\bbZ}H^i(G,k)$ for $p=2$ and $\oplus_{i\in2\bbZ}H^i(G,k)$ for $p$ odd respectively (and stress the awkwardness of this definition). Let $\stab(kG)$ be the tensor triangulated category of finitely generated (left) $kG$-modules modulo projective modules, where $\otimes=\otimes_k$. Let us denote by $\Pi(G):=\Proj(H^\smallbullet(G,k))$ and for any finitely generated $kG$-module, by $\sigma(M):=\Pi(G)_M$ the $\Pi$-support of $M$ as defined in~\cite[Def.\,4.1]{fripev}.
\end{Not}

\begin{Thm}
\label{fripev-thm}
The pair $\big(\,\Pi(G)\,,\,\sigma\big)$ is a classifying support data on~$\stab(kG)$.
\end{Thm}

\begin{proof}
See~\cite[Thm.\,5.3]{fripev} for finite group schemes or~\cite[Thm.\,3.4]{bcr} for finite groups\,; in fact $M\mapsto\Pi(G)_M$ is a support data by~\cite[Prop.\,4.2]{fripev} or~\cite[Prop.\,2.2]{bcr}.
\end{proof}

\begin{Cor}
\label{modular-cor}
We have a homeomorphism $f:\Pi(G)\isoto\Spc\big(\stab(kG)\big)$ with
$$
\qquad
f(x)=\SET{a\in\stab(kG)}{\,x\notin\Pi(G)_M}
\qquad\text{for all }x\in\Pi(G)\,.
$$
Moreover, for any finitely generated $kG$-module $M$, the closed subset $\Pi(G)_M\subset\Pi(G)$ corresponds via $f$ to $\supp(a)$ in $\Spc\big(\stab(kG)\big)$.
\end{Cor}

\begin{proof}
Theorems~\ref{exogene-thm}, \ref{fripev-thm} and~\ref{univ-thm}.
\end{proof}


\bigbreak
\section{The structure sheaf}
\label{ring-sect}
\medbreak


%
\begin{Def}
\label{ring-def}
For any open $U\subset\Spccat{K}$, let $Z:=\Spccat{K}\oursetminus U$ be its closed complement and let $\KZ{K}{Z}$ be the \tthsub\ of $\cat{K}$ supported on $Z$ (see~\ref{KZ-not}). We denote by $\calOcat{K}$ the sheafification of the following presheaf of rings\,:
$$
U\mapsto \End_{\cat{K}/\KZ{K}{Z}}(\unit_U)
$$
where the unit $\unit_U\in\cat{K}/\KZ{K}{Z}$ is the image of the unit $\unit$ of $\cat{K}$ via the localization $\cat{K}\to\cat{K}/\KZ{K}{Z}$. The restrictions homomorphisms in the above presheaf are given by localization in the obvious way. The sheaf of commutative rings $\calOcat{K}$ turns $\Spccat{K}$ into a ringed space, that we denote\,:
$$
\Speccat{K}:=\big(\Spccat{K}\,,\,\calOcat{K}\big)\,.
$$
\end{Def}

\begin{Rem}
\label{presheaf-rem}
This construction is inspired by the author's~\cite{bal}. There, we considered a presheaf of triangulated categories on $\Spccat{K}$ given by $U\mapsto\widetilde{\cat{K}/\KZ{K}{Z}}$, with $Z=\Spccat{K}\oursetminus U$ (for the idempotent completion $\iccat{K}$, see Rem.\,\ref{ic-tens-rem}). This was in turn inspired by Thomason's theorem~\cite[Thm.\,2.13]{bal}, which identifies the latter category with $\Dperf(U)$ when $\cat{K}=\Dperf(X)$. Note that the endomorphism ring of the unit is commutative, see \eg~\cite[Lem.\,9.6]{bal}, and that the idempotent completion is harmless for the definition of the structure sheaf $\calOcat{K}$, since $\cat{K}\hookrightarrow\iccat{K}$ is always a full embedding.

In~\cite{bal}, we established that this sheaf of rings recovers the structure sheaf $\calO_{X}$ when applied to $\cat{K}=\Dperf(X)$. We will not repeat this in the present context and we leave it as an easy exercise to the reader, considering the above comments. The computation of the ``right'' structure sheaf in the case of $\cat{K}=\stab(kG)$ comes more as a surprise because of the rather non-conceptual definition of $H^\smallbullet(G,k)$, see~\ref{group-not} and compare Benson~\cite[Vol.\,II, Rem.\ after 5.6.4, p.\,175]{ben}. Nevertheless, we have the following result.
\end{Rem}

\begin{Thm}
\label{examples-thm}
Via the homeomophisms of Corollaries~\ref{scheme-cor} and~\ref{modular-cor}, the structure sheaves also identify. That is, we have the following isomorphisms of schemes\,:
\begin{enumerate}
\item For $X$ a topologically noetherian scheme, $\Spec\big(\Dperf(X)\big)\simeq X$.
\smallbreak
\item For $G$ a finite group (scheme), $\Spec\big(\stab(kG)\big)\simeq \Proj\big(H^\smallbullet(G,k)\big)$.
\end{enumerate}
\end{Thm}

\begin{proof}
As already mentioned, Part~(a) follows as in~\cite{bal} from Thomason's theorem, see~\cite[Thm.\,2.13]{bal}. Part~(b) has been established recently by Friedlander and Pevtsova and is precisely the statement of~\cite[Thm.\,7.3]{fripev}.
\end{proof}

\begin{Rem}
\label{final-rem}
It is an open question to know when $\Speccat{K}$ is a scheme. For the moment, we can prove that $\Speccat{K}$ is always a \emph{locally} ringed space. Since we do not have applications of this fact yet, we do not include its proof here.
\end{Rem}


\bibliographystyle{abbrv}


\end{document}